\documentclass{article}

\usepackage{arxiv}

\usepackage[utf8]{inputenc} 
\usepackage[T1]{fontenc}    
\usepackage{hyperref}       
\usepackage{url}            
\usepackage{booktabs}       
\usepackage{amsfonts}       
\usepackage{microtype}      
\usepackage{graphicx}
\usepackage{natbib}

\usepackage{amssymb}
\usepackage{amsmath}

\newtheorem{theorem}{Theorem}

\newtheorem{lemma}{Lemma}
\newtheorem{corollary}{Corollary}

\newtheorem{proposition}{Proposition}
\newtheorem{remark}{Remark}
\newtheorem{example}{Example}

\newcommand{\e}{\text{\large{$\mathbf{j}$}}} 
\newenvironment{proof}[1][Proof]{\textit{\textbf{#1:}} }

\title{Main  signless Laplacian eigenvalues of \emph{quasi}-threshold graphs}

\author{\'{A}tila Jones \\
Instituto Federal Sudeste De Minas \\
Campus Juiz De Fora, Juiz De Fora, MG, Brazil.\\
	\texttt{atila.jones@ifsudestemg.edu.br} \\
	\And Vilmar Trevisan\\
	Instituto de Matem\'{a}tica e Estat\'{\i}stica\\
 Universidade Federal do Rio Grande do Sul  Porto Alegre, RS, Brazil.\\
	\texttt{trevisan@mat.ufrgs.br} \\
\And
Cybele T.M.Vinagre \\
	 Instituto de Matem\'atica e Estat\'{\i}stica \\
	 Universidade Federal Fluminense, Niter\'oi, RJ, Brazil \\
	\texttt{cybele\_vinagre@id.uff.br} \\
}

\renewcommand{\headeright}

\hypersetup{
pdftitle={Main Q-eigenvalues of quasi-threshold graphs},
pdfsubject={combinatorics},
pdfauthor={\'{Atila} Jones, Vilmar Trevisan, Cybele T.M. Vinagre},
pdfkeywords={Cograph,\emph{Quasi}-threshold graph,  Signless Laplacian matrix,  Main signless Laplacian eigenvalue}
}

\begin{document}
\maketitle

\begin{abstract}
	In this note, we present a structural description of certain connected cographs having $k \geq 2$ main signless Laplacian eigenvalues. This result allows us to characterize the cographs which are \emph{quasi}-threshold graphs with two main $\mathbf{Q}$-eigenvalues. In addition, we describe all the \emph{quasi}-threshold graphs belonging to the subclass of generalized core-satellite graphs with $k \geq 2$ main $\mathbf{Q}$-eigenvalues.
\end{abstract}


\section{Introduction and main results}
\label{sec:intro}

An eigenvalue $\lambda$ is a \emph{main eigenvalue} of a square matrix $\mathbf{M}$ (or a \emph{main $\mathbf{M}$-eigenvalue}) if the eigenspace $\mathcal{E}(\lambda)$ of $\lambda$ is non-orthogonal to the all ones vector $\e$, that is, if $\mathcal{E}(\lambda)$ contains  some eigenvector whose sum of entries is non-zero.

Main eigenvalues of graphs were introduced  in \citep*{cvet1970} for the case $\mathbf{M} =\mathbf{A}(G)$, the adjacency matrix of the graph $G$.  In \citep*{ChenHuang2013} and \cite{denghuang2013}, the case $\mathbf{M} =\mathbf{Q}(G)$, the signless Laplacian matrix of  $G$, was studied and it was shown that regular graphs were the only ones with exactly one main $\mathbf{Q}$-eigenvalue. The problem of characterizing graphs with $k \geq 2$ main $\mathbf{Q}$-eigenvalues was also addressed, in particular the case $k=2$ for trees and unicyclic graphs. In \citep*{ChenHuang2013}, the bicyclic graphs with  two main $\mathbf{Q}$-eigenvalues were characterized. In addition, trees, unicyclic and bicyclic graphs with three main $\mathbf{Q}$-eigenvalues, one of which being zero, were characterized in \citep*{Javarsineh2017} and \citep*{JAVARSINEH2017603}.

In \citep{VINAGRE202033}, the authors gave an exact characterization of the threshold graphs having  $k$ main $\mathbf{Q}$-eigenvalues, for each $k \geq 2$. The present paper is a step forward to generalize this characterization for \emph{quasi}-threshold graphs.

We may recall that threshold graphs form a subclass of cographs characterized as the class of $\{P_4, C_4, 2K_2\}$-free graphs, while \emph{quasi}-threshold graphs constitute the class of $\{C_4, P_4\}$-free graphs (we refer to the definitions below).

The main result of this paper is a structural characterization of connected \emph{quasi}-threshold graphs having $k\geq 2$ main $\mathbf{Q}$-eigenvalues. In addition, we determine precisely what are the connected \emph{quasi}-threshold graphs having two main $\mathbf{Q}$-eigenvalues. We also provide a full characterization of a subclass of \emph{quasi}-threshold graphs. More precisely, for $k=2,3,\ldots $, we determine what are the generalized core-satellite graphs having $k$ main $\mathbf{Q}$-eigenvalues. It is important to note that the graph characterizations shown in this article were obtained with the help of the \textit{Graph Filter} software \citep{graphfilter}.

The paper is organized as follows. In the next section, we provide the definitions and the necessary background  for presenting the results.  We then start, in Section~\ref{sec:bound}, studying the number $k$ of main $\mathbf{Q}$-eigenvalues of a connected  cograph $G$, presenting an upper bound on $k$ based on the cotree representation of $G$. This allows us to compute the spectra of some particular \emph{quasi}-threshold graphs, a necessary step to  the results of the next section. In fact, in Section~\ref{sec:char}, we give a structural characterization of the \emph{quasi}-threshold graphs having $k \geq 2$ main $\mathbf{Q}$-eigenvalues. We show they have the form $G=K_c \oplus H$ for some integer $c\geq 1$ and a disconnected \emph{quasi}-threshold graph $H$. Additionally, we prove that if $\overline{H}$ is bipartite, then $k=2$. In turn,  these results lead to a precise determination of the connected \emph{quasi}-threshold graphs having two main $\mathbf{Q}$-eigenvalues given in Theorem~\ref{theo:chordal2main}. Section~\ref{sec:sat} is devoted to the study of main $\mathbf{Q}$-eigenvalues of generalized core-satellite graphs, an important subclass of the class of \emph{quasi}-threshold graphs, which contains, for example,  the windmill graphs. For this class, we characterize which graphs have $k$ main $\mathbf{Q}$-eigenvalues, for $k=2, 3, \ldots$.

\section{Preliminaries}\label{sec:prel}

Let $G =(V, E)$ be a simple  graph. We denote by $d_G(u)=d(u)$ the degree of vertex $u \in V$. As usual, we say that $|E| =m$ is the \emph{size} of $G$ and $|V| =n$ is the \emph{order} of $G$.
The \emph{neighborhood of a vertex } $v\in V$ is the set $N_G(v) =\{w \in V; \{v, w\} \in E \}$ of neighbors of $v$ in $G$ and  the \emph{closed neighborhood}  of $v$ is  $N_G[v] =N_G(v)\cup \{v\}$.

The complete graph on $n$ vertices is denoted by $K_n$. Its complement $\overline{K_n}$ is the graph without edges. For an integer $p \geq 1$ and a graph $G$, $pG$ denotes the disjoint union of $p$ copies of $G$. For any $S\subset V$, the subgraph of $G$ induced by $S$ is denoted $G[S]$. If $G[S]$ is a complete subgraph then $S$ is a \emph{clique in} $G$. A vertex $v$ is said to be \emph{simplicial} (respectively, \emph{universal})  in $G$ when $N_G(v)$ is a clique (respectively,  $N_G[v] =V$).

\subsection{Chordal graphs}

A \emph{chordal graph} is a graph in which every cycle of length four or greater has a \emph{chord}, that is, an edge between non-consecutive vertices of the cycle. Basic concepts about chordal graphs 
can be found in \citep{blair1993introd}.
The notions which are relevant to this note follow.
A subset $S\subset V$ is a \emph{separator} of a non-complete graph $G$ if its removal leaves the remaining graph disconnected.
The set $S$ is a \emph{minimal separator} of $G$ if $S$ is a separator and no proper subset of $S$ separates the graph. The minimum cardinality of a separator of $G$ corresponds to $\kappa(G)$, the \emph{vertex connectivity} of the graph.
Let $u, v$ be two non-adjacent vertices in $V$. A subset $S\subset V$ is a $uv$-\emph{separator vertex} (or a \emph{separator for} $u$ \emph{and} $v$)  if the removal of $S$ from the graph separates $u$ and $v$ into distinct connected components. If no proper subset of $S$ is a $uv$-separator then $S$ is a \emph{minimal $uv$-separator}. When the pair of vertices remains unspecified, we refer to $S$ as a \emph{minimal vertex separator} (\textit{mvs}). A minimal separator is always a minimal vertex separator but the converse is not true.

\begin{lemma}\label{lem:simpl} 
Suppose that a  connected regular graph $G$ has a simplicial vertex. Then $G$ is a complete graph.
\end{lemma}

\begin{proof}
Let $u$ be a simplicial vertex of the  graph $G$ and suppose that $S=N_G[u]$ is a $k$-clique; then $d(u)=k-1$. Suppose there is a vertex $v$ in $G\setminus S$. Since $G$ is a connected graph, there exists a vertex $w \in S$, $w \not= u$, and there exists a path  $vv_1\ldots v_rw$ in $G$. Then $d(w) \geq (k-1)+1=k$, in contradiction with the regularity of  $G$. Therefore, $G\setminus S$ is the empty set and $G$ is a complete graph.
\end{proof}
\begin{corollary}\label{cor:cor1} A connected regular chordal graph is a complete graph.
\end{corollary}
\begin{proof}
    It is well known that a graph $G$ is chordal if and only if every induced subgraph of $G$ has a simplicial vertex. So, by Lemma \ref{lem:simpl}, the result follows.
\end{proof}

\subsection{Cographs}

The class of \emph{complement reducible graphs}, known as \emph{co\-graphs}, is the hereditary class of $P_4$-free graphs, denoted by $\mathcal{F}orb(P_4)$.
This class appeared naturally in various contexts and has been rediscovered several times, but  the seminal paper   \citep*{Corneil1981} is considered the origin of the study of cographs.
There are several equivalent ways in which the elements of this class may be defined. One of the definitions uses disjoint unions and complements recursively, as the authors of \citep{Corneil1981}.
A similar definition may be obtained through unions and joins. Consider the class $\mathcal{D}$ defined recursively as follows:
\begin{itemize}
\item[(${\alpha}$)]~$K_1 \in \mathcal{D}$;

\item[(${\beta}$)]~If $G_1$ and $G_2$ are vertex disjoint graphs  in $\mathcal{D}$, then $G_1 \cup G_2 \in \mathcal{D}$;

\item[(${\gamma}$)]~If $G_1$ and $G_2$ lie in $\mathcal{D}$, then its join $G_1 \oplus G_2 \in \mathcal{D}$.
\end{itemize}

Using this definition, it is possible to characterize bipartite cographs (important for our purposes in this note) as the following known result shows.

\begin{lemma}\citep[Lemma~14]{DEMANGE2005} \label{lem:cogbip}
  If $G$ is a connected and bipartite cograph then \  $G=\overline{K_a}\oplus \overline{K_b}$ \  for positive integers $a, b $.
\end{lemma}

We recall that \emph{threshold graphs} constitute  a subclass of cographs characterized as the hereditary class of $\{P_4, C_4, 2K_2\}$-free graphs. Since, by definition,  chordal graphs are $C_4$-free graphs, the class of chordal cographs is the hereditary class of $\{C_4, P_4\}$-free graphs whose elements are often called \emph{quasi-threshold graphs} and constitute  an important class of graphs (see, for example, \citep{YanChen1996}).

\begin{remark}\label{rem:cognonchord} 
  For fixed positive integers $a_1, a_2, b_1,b_2$, we may note that  $(K_{a_1}\cup K_{b_1})\oplus (K_{a_2}\cup K_{b_2})$ is a non-chordal cograph. To see this, it suffices to choose one vertex at each different copy of the complete graph to obtain an  induced $C_4$. In particular, if $c \geq 2$, then the graph of the form $\overline{K_c} \oplus  (K_{a_2}\cup K_{b_2})$ is a non-chordal cograph.
  On the other hand,  for $t \geq 1$ and $a_i \geq 1$, $1 \leq i\leq t$, the graph $t K_{a_i}$ is a disconnected chordal cograph.


\end{remark}

\subsection{Signless Laplacian matrix and their main eigenvalues}

For a simple graph  $G$ with vertices $v_1, \ldots, v_n$\,,  its adjacency matrix $\mathbf{A}(G) =[a_{ij}]$ is the  square matrix of order $n$ where $a_{ij}=1$ if $v_i$ and $v_j$ are adjacent and 0, otherwise. The \emph{signless Laplacian matrix} $\mathbf{Q}(G)$ is defined as $\mathbf{Q}(G) =\mathbf{D}(G) +\mathbf{A}(G)$, where $\mathbf{D}(G)$ is the diagonal matrix whose entries are the degrees of the vertices. The matrix $\mathbf{Q}=\mathbf{Q}(G)$ is symmetric and positive semidefinite and then it has $n$ real non-negative eigenvalues, which will be called $\mathbf{Q}$-\emph{eigenvalues of } $G$. We denote the multiset of these eigenvalues by $\sigma(G)=\{q^{[r_1]}_1, q^{[r_2]}_2, ..., q^{[r_p]}_p\}$, where
and $q^{[r]}$ means that the eigenvalue $q$ has multiplicity $r\geq 1$.

The following known property of the matrix $\mathbf{Q}(G)$ of a bipartite graph  $G$ is stated here for future reference.

\begin{lemma}\citep{CVETKOVIC2007}\label{lem:zeroasQeigenval} 
  For any connected graph  $G$\,, the least $\mathbf{Q}$-eigenvalue of  $G$ is 0 if and only if  $G$ is bipartite. In this case, 0 is a simple eigenvalue. In general, the multiplicity of 0 as a $\mathbf{Q}$-eigenvalue of $G$ is the number of bipartite components of  $G$.
\end{lemma}

An eigenvector $\mathbf{v}$ associated with a  $\mathbf{Q}$-eigenvalue $q$ of $G$   is a \emph{main eigenvector} of $G$ if  $\mathbf{v}^{\top}\e \neq 0$. In this case, $q$ is said to be a \emph{main $\mathbf{Q}$-eigenvalue} of $G$. From the Perron-Frobenius' Theorem for non-negative irreducible real symmetric matrices,  if  $G$ is a connected graph on $ n  \geq 2$ vertices then its largest $\mathbf{Q}$-eigenvalue is simple and main,  since it has an associated eigenvector whose entries are all positive.
The main $\mathbf{Q}$-eigenvalues of a graph may be considered as pairwise distinct eigenvalues, since it is possible to find an orthonormal basis for each one of the associated eigenspaces containing only a single eigenvector non-orthogonal to $\e$. This fact,   whose proof can be found in \citep[Remark~3]{Cardoso2012}, will be taken into account throughout the paper.

The  facts about main $\mathbf{Q}$-eigenvalues of a graph listed below will be applied in the remaining of the text.

\begin{lemma} \citep*{denghuang2013}, \citep*{ChenHuang2013} \label{lem:reg}
 A graph $G$ has exactly one main signless Laplacian eigenvalue if and only if it is regular.

\end{lemma}

\begin{lemma}\citep{VINAGRE202033} \label{lem:maincompl} 
  A graph  $G$ and its complement $\overline{G}$ have the same number of main $\mathbf{Q}$-eigenvalues.
\end{lemma}

The $\mathbf{Q}$-spectrum of the union of two graphs is the union of the $\mathbf{Q}$-spectra of the original graphs  \citep[see][]{Cvetkovic2009}.   It follows that the main $\mathbf{Q}$-eigenvalues  of the union of two graphs are obtained from the union (here considered as a set, not a multiset) of the main-$\mathbf{Q}$ eigenvalues of  each one of the graphs.
The following  two results are consequences of these facts.

\begin{lemma}\citep{VINAGRE202033} \label{lem:non-zeroismain}
  Let $H$ and $\overline{K_p}$\,, for $p\geq 1$ an integer,   be graphs with disjoint vertex sets. A real number $q\neq 0$ is a main $\mathbf{Q}$-eigenvalue of $H$ with main eigenvector $\mathbf{v}$ if and only if $q$ is a main $\mathbf{Q}$-eigenvalue of the union $G =\overline{K_{p}}\ \cup $ H with associated main eigenvector
  $\left[\begin{array}{c}                                                              \mathbf{v} \\                                                                        \mathbf{0}_p\\                                                                       \end{array}\right]$,
  where $\mathbf{0}_p$ stands for the null vector on $p$ coordinates.

\end{lemma}

The following result was proven for a connected $H$, but holds for a non-connected graph as well.
\begin{lemma}\citep{VINAGRE202033} \label{lem:zeromainofcompl}
  Any graph $G =\overline{K_p}\ \cup H$,  where  $p  \geq 1$ is an integer and $H$ is a graph, has 0 as a main $\mathbf{Q}$-eigenvalue.

\end{lemma}

\section{A bound on the number of main $\mathbf{Q}$-eigenvalues of a cograph}\label{sec:bound}

In this section, we use a cograph representation to provide a bound on the number of main $\mathbf{Q}$-eigenvalues of cographs. Before doing so, we use  Lemma~\ref{lem:zeromainofcompl} to prove a result that predicts the number of main $\mathbf{Q}$-eigenvalues of a particular cograph.
\begin{theorem}\label{tent1}
 Let $G$ be cograph with $k \geq 2$ main $\mathbf{Q}$-eigenvalues. Then
  \begin{enumerate}
    \item[$(a)$] For any integer $c\geq 1$, the cograph $K_c \oplus G$ has $k$ or  $k+1$ main $\mathbf{Q}$-eigenvalues;
    \item[$(b)$] For each integer $c\geq 1$, the cograph $K_c \oplus G$ has  $k+1$ main $\mathbf{Q}$-eigenvalues if and only if $\overline{G}$ is not  bipartite.
    \end{enumerate}
\end{theorem}

\begin{proof}
  Let $G$ be a cograph with $k \geq 2 $ main $\mathbf{Q}$-eigenvalues and let $c\geq 1$ be an integer.  By Lemma \ref{lem:maincompl}, we know  that $\overline{G}$ has $k$ main $\mathbf{Q}$-eigenvalues. By invoking Lemma~\ref{lem:zeromainofcompl}, we notice that $\overline{K_c} \cup \overline{G}$ has 0 as a main $\mathbf{Q}$-eigenvalue. Regarding item (a), we need to consider whether 0 is a  main $\mathbf{Q}$-eigenvalue of $\overline{G}$ or not. Since we are considering that the main $\mathbf{Q}$-eigenvalues are all distinct,  in the first case we have that  $\overline{K_c} \cup \overline{G}$ remains with $k$ main $\mathbf{Q}$-eigenvalues and, in the second, that the cograph $\overline{K_c} \cup \overline{G}$ has $k+1$ main $\mathbf{Q}$-eigenvalues. Applying again Lemma \ref{lem:maincompl}, the assertion of item (a) is proved. Taking into account Lemma \ref{lem:zeroasQeigenval},  the previous argumentation   leads us  to conclude also that the item (b) is true.
\end{proof}

We may represent every cograph on vertex set $\{1,\ldots, n\}$ as a rooted tree, named  \emph{cotree}, whose nodes  consist of $n$ leaves labeled $1$ through $n$ and internal vertices that carry either the label ``$\cup$'' for union or ``$\oplus$'' for join. Given such a tree $T$, we can easily construct the corresponding cograph $G_T$. The construction can be made unique (see, for example \citep{livro}) as follows. We say that a cotree is in its  \emph{normalized form} if every internal node has at least two children and has a label that differs from the label of its parent. In other words, the children of nodes labeled $\cup$ are leaves or nodes labeled $\oplus$, while the children of nodes labeled $\oplus$ are leaves or nodes labeled $\cup$. In this note, we assume that the cotree representation of a cograph is always in normalized form.

Consider the cotree $T_G$ of a cograph $G$.
We group  the leaves of the cotree that represent the vertices of $G$ as follows.
Let $r$ be the number of distinct sets of siblings, which we call the \emph{width} of the cotree. We call each of these sets of  sibling  vertices  a \emph{bag} $B_i$ and denote the respective cardinality as $t_i=|B_i|$, \ $1 \leq i \leq r$. When the siblings are \emph{co-duplicates} (meaning the parent is $\oplus$), we call it a $J$-\emph{bag}, and when the siblings are \emph{duplicates} (meaning the parent is $\cup$), we call it  a $U$-\emph{bag}. If the cotree $T_G$ has $r$ bags $B_i$ whose cardinalities are $t_i$, we can represent $T_G$ as the sequence $$\left\{B_i^{t_i}\right\}_{i=1}^r,$$
which we call \emph{bag representation} of $T_G$.

We state the following result that uses this representation of cographs and is useful for our purposes.

\begin{proposition}\label{prop:nonmaincogr} Let $G$ be a cograph with cotree $T_G$ having $r_1$ $J$-bags $B_i^{t_i}$\,, $i=1,\ldots, r_1$\,, and $r_2=r-r_1$ $U$-bags $B_j^{t_j}$\,, $j=r_1+1,\ldots, r$. Let  $p_i$ be the degree of each of the vertices of $G$ in the bag $B_i, ~i=1,\ldots ,r$. Then
\begin{enumerate}
  \item[$(a)$] $q=p_i-1$ and $q=p_j$ are eigenvalues of $\mathbf{Q}(G)$ with multiplicities at least $t_i-1$, for $i=1,\ldots, r_1$\,, and $t_j-1$, for $j=r_1+1,\ldots,r$, respectively.
  \item[$(b)$] All of them are non-main $\mathbf{Q}$-eigenvalues of $G$.
  \end{enumerate}
\end{proposition}

\begin{proof} We notice that item (a) is given in \citep[Proposition 2(c)]{JONES2023120}. In any case, the proof of item (b) we give now, also proves item (a).
For a cograph $G$ as above, where each vertex within the bag $B_i$ has
degree equal to $p_i$, we consider the matrix $\mathbf{Q}(G)$.  For each pair of  vertices $u, v$ belonging to the same  $J$-bag  (respectively, $U$-bag) $B_i$,  we consider the vector whose entries are all equal to 0 except those corresponding to $u$ and $v$, which are taken to be 1 and -1, respectively. We obtain a $\mathbf{Q}(G)$-eigenvector associated with eigenvalue $p_i-1$ (respectively, $p_i$). This reasoning allows one to conclude that, for each $i=1,\ldots, r_1$,  the integer $p_i-1$ is a non-main $\mathbf{Q}$-eigenvalue of the graph $G$ with multiplicity at least $t_i-1$ and for each $j=r_1+1,\ldots,r$, $p_j$ is a  non-main $\mathbf{Q}(G)$-eigenvalue  with multiplicity at least $t_j-1$.
\end{proof}

We observe that, from Proposition~\ref{prop:nonmaincogr}, we obtain \ $n-r$ \ $\mathbf{Q}$-eigenvalues of $G$ that are non-main.  Further, as  proved in \citep{JONES2023120}, the remaining $r$ eigenvalues are those of the $r\times r  $ \emph{condensed signless Laplacian matrix}, which is the symmetric matrix defined as  $\mathbf{C}(G)=[c_{ij}]$, where
 \[c_{ij} =\left\{\begin{matrix}
p_{i}+(t_i-1), & \mbox{ if } & i=j \mbox{ and } 1 \leq i \leq r_1; \\
p_i, & \mbox{ if } & i=j \mbox{ and } r_1 +1 \leq i \leq r;\\
\sqrt{t_it_j}z_{ij} &\mbox{ if } & i\neq j \mbox{ and } 1 \leq i,j \leq r,
\end{matrix}\right. \,,\]
where $z_{ij}=1 \mbox{ or }  0$, depending on whether or not  the vertices of bag $B_i$ are adjacent to the vertices of $B_j$. Thus, the matrix $ \mathbf{C}(G)$ contains the main $\mathbf{Q}$-eigenvalues of the graph $G$ and we obtain the following result as a consequence.

\begin{theorem}\label{the:boundmaincog}
If $G$ is a cograph with bag representation $$\left\{B_i^{t_i}\right\}_{i=1}^r,$$
  then $G$ has at most $r$ main $\mathbf{Q}$-eigenvalues.
\end{theorem}

Hence, the number of main $\mathbf{Q}$-eigenvalues of a cograph does not exceed the width of its cotree.  We notice that, in general, this bound can be far from optimal.  For example, as we will see in Lemma \ref{lem:mainUnion}, the cograph $K_c\oplus(tK_a)$, for integers $a, c\geq 1$ and  $t \geq 2$ has two main $\mathbf{Q}$-eigenvalues while its width, which is $t+1$, can be arbitrarily large.

In the sequence, we apply the  facts mentioned in the proof of Proposition \ref{prop:nonmaincogr} to identify the main $\mathbf{Q}$-eigenvalues  of a connected bipartite cograph that, according to Lemma~\ref{lem:cogbip}, has the form $\overline{K_a}\oplus \overline{K_b}$ for some integers $1\leq a\leq b$.

Firstly, we may recall from \citep{GroneMerris1990LaplacianSpectrum} that the $\mathbf{Q}$-spectrum of the complete graph on $n$ vertices is
\[\sigma(K_n)=\{(2n-2)^*, (n-2)^{[(n-1)]}\}\,, \]
where  $q^*$ indicates that $q$ is a main $\mathbf{Q}$-eigenvalue of the graph.

For bipartite graphs, as  the $\mathbf{Q}$-characteristic polynomial and the $\mathbf{L}$-characteristic polynomial are equal (definition of Laplacian matrix in the next section), it is known (see, for example, \citep{CVETKOVIC2007} and \citep{GroneMerris1990LaplacianSpectrum}) that
$\sigma(\overline{K_a}\oplus \overline{K_b})=\{(a+b)\,,\  b^{[a-1]}\,,\ a^{[b-1]}\,, 0\}$, for $a,b$ non-negative integers.
 If $a,b\geq 1$ are distinct integers, we claim that the main $\mathbf{Q}$-eigenvalues of $ G= \overline{K_a}\oplus \overline{K_b}$ are $a+b$ and 0. Indeed, from the proof of Proposition \ref{prop:nonmaincogr}, the condensed signless Laplacian matrix of this graph is
  \[\mathbf{C}(\overline{K_a}\oplus \overline{K_b})= \left[
      \begin{array}{cc}
        b & \sqrt{ab} \\
        \sqrt{ab} & a \\
      \end{array}
    \right],
  \]
whose eigenvalues are exactly $a+b$ and $0$. From Theorem~\ref{the:boundmaincog}, the number $k$ of main $\mathbf{Q}$-eigenvalues is at most 2. But since $G$ is not a regular graph, we have that $k =2$, by Lemma~\ref{lem:reg}. In case $a=b$, the graph is regular and there is a single main $\mathbf{Q}$-eigenvalue, which is $2a$. For future reference, we state the following result.

\begin{lemma}\label{lem:spectbipcomp}
If $a,b$ are positive integers with $a \neq b$ then
\[\sigma(\overline{K_a}\oplus \overline{K_b})=\{(a+b)^*\,,\  b^{[a-1]}\,,\ a^{[b-1]}\,, 0^*\}\,.\]
 In the particular case when $a=b$, the graph is regular and
\[\sigma(\overline{K_a}\oplus \overline{K_a})=\{(2a)^*\,,\  a^{[2a-2]}\,, 0\}\,.\]
\end{lemma}

From this result, we obtain the number of main $\mathbf{Q}$-eigenvalues of two  types of cographs which we will see, in Proposition \ref{prop:qt-con} below, that are \emph{quasi}-threshold graphs.

\begin{lemma}\label{lem:mainUnion} Consider the cographs
  \begin{align*}
  G_1& = K_c \oplus (K_a \cup K_b), \  \mbox{ for integers }  a, b, c \geq 1  \mbox{ with } a\not = b  \  \ \mbox{ and } \\
     G_2&=K_c\oplus (t K_a), \mbox{ for integers } \mbox{ for }  t \geq 2, \  a, c \geq 1.
     \end{align*}
     The number of main $\mathbf{Q}$-eigenvalues of $G_1$ and  $G_2$ is 2.
\end{lemma}
\begin{proof}
In fact,  the signless Laplacian condensed matrix  $\mathbf{C}(G_1)$ is the $3\times 3$ matrix
   \[\left[
       \begin{array}{ccc}
         a+b +2(c-1) & \sqrt{ca} & \sqrt{cb} \\
         \sqrt{ca} & 2(a-1)+c & 0 \\
         \sqrt{cb} & 0 & 2(b-1)+c \\
       \end{array}
     \right]\,,
   \]
  for which the eigenvalues are $a+b+c-2$ and the two roots of the equation
  \[q^2- (3c+2b+2a-4)q+2c^2 +(2b+2a-6)c+(4a-4)b-4a+4=0.\]
From Proposition  \ref{prop:nonmaincogr}, we see  that $q=a+b+c-2$ is one of the non-main $\mathbf{Q}(G_1)$-eigenvalues  and that its  multiplicity is at least  $c-1$,  since it corresponds to the degrees of vertices  of the $J$-bag with $c$ vertices. We claim that  $q=a+b+c-2$ actually has multi\-plicity equal to $c$ as a non-main eigenvalue of $\mathbf{Q}(G_1)$.  Indeed, it is not difficult to see that  $\mathbf{C}(G)\mathbf{v}=(a+b+c-2)\mathbf{v}$, for
 $\mathbf{v}=\left[                                                                \begin{array}{c}
 1 \\
  {c}/{(b-a)} \\                                                                 -{c}/{(b-a)} \\                                                                \end{array}
   \right]  $.
    Let us denote by $\mathbf{J}_{r\times s}$ (respectively, $\mathbf{O}_{r\times s}$) the all ones matrix (respectively, the null matrix) of order $r\times s$, by $\e_s$ the all ones vector of length $s$ and by  $\mathbf{I}_{r}$ the identity matrix of order $r$. Then, we can describe $\mathbf{Q}(G_1)$ as
  \[\mathbf{Q}(G_1)= \left[
     \begin{array}{ccc}
       \mathbf{J}_{c\times c} + (a+b+c-2)\mathbf{I}_{c}& \mathbf{J}_{c\times a} & \mathbf{J}_{c\times b} \\
       \mathbf{J}_{a\times c} & \mathbf{J}_{a\times a} +(a+c-2)\mathbf{I}_{a} & \mathbf{O}_{a\times b} \\
       \mathbf{J}_{b\times c} & \mathbf{O}_{b\times a} & \mathbf{J}_{b\times b} +(b+c-2)\mathbf{I}_{b} \\
     \end{array}
   \right]
  \]
and verify that     $\mathbf{w}=\left[                                                                \begin{array}{c}
 \e_c \\
  {c}/{(b-a)}\e_a \\                                                                 -{c}/{(b-a)}\e_b \\                                                                \end{array}
   \right]  $ is an eigenvector of $\mathbf{Q}(G_1)$ associated with the eigenvalue $q=a+b+c-2$, which forms a linearly independent set when grouped with the other eigenvalues associated with $q$ provided by Proposition \ref{prop:nonmaincogr}. Furthermore,
     the coordinates of  $\mathbf{w}$   add up to zero. This completes our claim. Since $G_1$ is not regular, the two remaining eigenvalues of $\mathbf{C}(G_1)$ must be  main $\mathbf{Q}$-eigenvalues of $G_1$. The assertion is proved for $G_1$.

  For $G_2$, by Lemma~\ref{lem:maincompl}, it suffices to compute the number of main $\mathbf{Q}$-eigenvalues of $\overline{G_2}= \overline{K_c} \cup \left(\oplus_{i=1}^t \overline{K_a}\right)$.
  Now,  $\overline{K_c}$ has 0 as the only main $\mathbf{Q}$-eigenvalue. As for $\left(\oplus_{i=1}^t \overline{K_a}\right)$, it is a regular graph having also a unique main $\mathbf{Q}$-eigenvalue. The degree of regularity is $(t-1)a$,  which is non-zero, since $t>1$, hence its main $\mathbf{Q}$-eigenvalue is non-zero, proving the result for $G_2$.
\end{proof}

\section{\emph{Quasi}-threshold graphs with $k$ main $\mathbf{Q}$-eigenvalues}\label{sec:char}
The \emph{Laplacian matrix} of a graph $G$ is defined as $\mathbf{L}(G) =\mathbf{D}(G)-\mathbf{A}(G)$. As $\mathbf{L}(G)$ is symmetric and positive semidefinite,  its  eigenvalues are  non-negative real numbers.   Fiedler  showed  in \citep{Fiedler1973} that $G$ is a connected graph if and only if  the second smallest eigenvalue of $\mathbf{L}(G)$ is  positive; for  this reason, this  eigenvalue is called the \emph{algebraic connectivity} of $G$ and  denoted by $a(G)$. Moreover, Fiedler proved that for $G\neq K_n$, $a(G) \leq  \kappa(G)$, where $\kappa(G)$ is  the vertex connectivity of $G$.

Graphs for which $a(G) = \kappa(G)$ were characterized  
in \citep{Kirkland2002} taking into account  their structures. It is known that $\kappa(G) =a(G)$ for cographs. A proof of this result can be found in \citep{ClaRenMarJoi} and we state it  here for easy reference.

\begin{lemma}\citep{ClaRenMarJoi}\label{lem:algconncog} Let $G$ be a cograph. Then $$\kappa(G) = a(G).$$
\end{lemma}

The particular case of chordal graphs for which $a(G) = \kappa(G)$ was treated in \cite{ABREU202168}
and we state the result here as the following lemma.

\begin{lemma} \citep[Theorem~2]{ABREU202168} \label{lem:abreuchordal}
Let  $G$ be a non-complete connected chordal graph. Then $\kappa(G) =a(G)$ if and only if there is a minimal separator of $G$ such that all its elements are universal vertices.
\end{lemma}

The next result gives a structural characterization of \emph{quasi}-threshold graphs that are connected.
\begin{proposition}\label{prop:qt-con} Let $G$ be a connected  \emph{quasi}-threshold. If $G$ is a non-complete graph, then $G=K_c \oplus H$, for some integer $c\geq 1$ and a disconnected \emph{quasi}-threshold graph $H$.
\end{proposition}

\begin{proof}
  Let $G$ be a  cograph which is connected and chordal.  Since it is not complete, it follows, by Lemma \ref{lem:algconncog}, that $a(G)=\kappa(G)$ and thus,  following  Lemma \ref{lem:abreuchordal},  $G$ admits a minimal vertex separator $S$ in which all the vertices are universal, implying that $S$
  is a clique. Let $|S|=c$, $c \geq 1$, be the cardinality of $S$. We have $S=K_c$, and because all vertices of $S$ are universal, we can describe $G$ as $G=K_c\oplus H$, \ where $H$ is a disconnected \emph{quasi}-threshold graph, proving the result.
\end{proof}

We now give information about the structure of a connected \emph{quasi}-threshold graph $G$ having $k$ main $\mathbf{Q}$-eigenvalues. Clearly, Corollary~\ref{cor:cor1} and Lemma~\ref{lem:reg} imply that $k=1$ if and only if  $G $ is a complete graph.

\begin{theorem}\label{teo:qt-estrut}
For a connected \emph{quasi}-threshold graph $G$, let $G=K_c \oplus H$ be its decomposition given in Proposition~\ref{prop:qt-con}. If $k\geq 2$ is the number of main $\mathbf{Q}$-eigenvalues of $G$ then, exactly one of the following situations occurs:
\begin{enumerate}
\item [$(a)$]  $\overline{H}$ is non-bipartite, that is,  $\overline{H}$ has $k-1$ non-zero main $\mathbf{Q}$-eigenvalues \  \ or
\item [$(b)$] $\overline{H}$ is bipartite, and in this case, $k=2$.
\end{enumerate}
\end{theorem}

\begin{proof}
     We notice that $\overline{H}$ is connected, since $H$ is disconnected. By Lemma \ref{lem:zeromainofcompl}, we know that $\overline{G}=\overline{K_c}\cup \overline{H}$ has 0 as a main $\mathbf{Q}$-eigenvalue. Also, Lemma \ref{lem:non-zeroismain} asserts  that a real $q \neq 0$ is a main $\mathbf{Q}$-eigenvalue of $\overline{G}$ if and only if $q$ is a main $\mathbf{Q}$-eigenvalue of $\overline{H}$.  We may recall that $\overline{G}$ also has $k$ main $\mathbf{Q}$-eigenvalues, by Lemma \ref{lem:maincompl}. Thus, in accordance to all these facts, we need to consider the only two possible cases: (i) $\overline{H}$ has $k-1$ non-zero main $\mathbf{Q}$-eigenvalues; and (ii) $\overline{H}$ has $k$ main $\mathbf{Q}$-eigenvalues, one of which is 0.
 In case (i), we see that $\overline{H}$ can not be bipartite by Lemma~\ref{lem:zeroasQeigenval} and, moreover, that all its $k-1$ main $\mathbf{Q}$-eigenvalues are non-zero. This proves item (a).
As for case (ii), it follows, from Lemmas \ref{lem:zeroasQeigenval} and \ref{lem:abreuchordal}, that $\overline{H}$ is a connected bipartite chordal cograph, which is of the  form $\overline{K_a} \oplus \overline{K_b}$, for some positive integers $a,b$, by Lemma~\ref{lem:cogbip}. Moreover, Lemma~\ref{lem:spectbipcomp} ensures that $\overline{H}$ has exactly one non-zero main $\mathbf{Q}$-eigenvalue. Hence $G=K_c \oplus (K_a \cup K_b)$ for some integers $a,b,c \geq 1$, which has two main $\mathbf{Q}$-eigenvalues, proving item (b).
\end{proof}

If the number of main $\mathbf{Q}$-eigenvalues of a connected chordal cograph is 2, we can be more precise and determine exactly which graph it is.

\begin{theorem}\label{theo:chordal2main}
   $G$ is a connected \emph{quasi}-threshold graph with two main \- $\mathbf{Q}$-eigenvalues if and only if $G$ assumes one of the following expressions:
   \begin{enumerate}
     \item[$(a)$] $G= K_c \oplus (K_a \cup K_b)$ for integers $a, b, c\geq 1$ and $a \neq b$;
     \item[$(b)$]  $G=K_c\oplus (t K_a)$, for integers $t \geq 2$, $a, c \geq 1$.
   \end{enumerate}
\end{theorem}

\begin{proof}
  Let $G$ be a cograph which is connected and chordal and has two main $\mathbf{Q}$-eigen\-values.
  By repeating the argument in the proof of Theorem~\ref{teo:qt-estrut}, we conclude that $G=K_c\oplus H$, for some integer $c\geq 1$ and $H$ a disconnected \emph{quasi}-threshold graph. Moreover, the only possible cases are: (i) $\overline{H}$ has two  main $\mathbf{Q}$-eigenvalues, one of which is zero; and (ii) $\overline{H}$ has only one main $\mathbf{Q}$-eigenvalue.
Concerning case (i), it follows from Lemma \ref{lem:zeroasQeigenval} that $\overline{H}$ is a connected bipartite  cograph, and it is of the form $\overline{H}=\overline{K_a} \oplus \overline{K_b}$, by Lemma~\ref{lem:cogbip}, for some positive integers $a\not =b$ (if $a=b$ then  $\overline{H}$ would be regular and consequently with a single main $\mathbf{Q}$-eigen\-value). Hence $G=K_c \oplus (K_a \cup K_b)$ for some integers $a,b,c \geq 1$, $a\not =b$.
In case (ii), the assumption on $\overline{H}$ implies that $H$ is a regular graph, by Lemmas~\ref{lem:reg} and \ref{lem:maincompl}. So, $G=K_c\oplus H$, where $H$ is a cograph which is disconnected,  regular and chordal. In view of Corollary~\ref{cor:cor1} and Remark~\ref{rem:cognonchord}, $H$ is the union of at least two copies of the same complete graph $K_a$, where $a \geq 1$ is an integer (a unique copy of a complete graph would imply $G$ to be  complete, a contradiction). So $G= K_c \oplus (t K_a)$, for integers $t \geq 2$ and $a \geq 1$.
Now, combining these two cases, we have that $G$ has one of the forms of items (a) or (b).  Reciprocally, if $G$ has the forms of item (a) or item (b),  Lemma~\ref{lem:spectbipcomp} asserts that $G$ has two main  $\mathbf{Q}$-eigenvalues.
  \end{proof}

\begin{remark} \label{generthreshold}
  In  Theorem \ref{theo:chordal2main}, for an integer $c \geq 2$, the \emph{quasi}-threshold graphs having the form  $K_c\oplus (K_a \cup K_b)$, with  $a=1 $ or $b=1$,  together with those of the form $K_c\oplus \left( ~ t K_a\right)$, with $t \geq 2$ and $a=1$, are exactly the connected threshold graphs with two main  $\mathbf{Q}$-eigenvalues whose characterization was given in \cite{VINAGRE202033}. In the latter case, we have the graphs $G=K_c\oplus  \overline{K_t}$, the so called \emph{complete split graphs}.
\end{remark}


\section{Main $\mathbf{Q}$-eigenvalues of generalized core-satellite graphs}\label{sec:sat}
In \citep{ESTRADA2017}, the graphs composed by a central clique (the \emph{core}) connected to several other cliques of same size (the \emph{satellites}) were called \emph{core-satellite graphs};  they generalize both the complete split graphs and the windmill graphs. The authors of \citep{ESTRADA2017} also considered the generalization of core-satellite graphs to the case where the satellites cliques are not restricted to having all the same order; they called these graphs \emph{generalized core-satellite}. It is not difficult to see that these graphs are also \emph{quasi}-threshold graphs. This fact was pointed out in \citep{YanChen1996}.

Then the general form of a generalized core-satellite graph is $$G=K_{n_0}\oplus (a_1K_{n_1} \cup a_2K_{n_2}\cup \ldots \cup a_p K_{n_p})=K_{n_0}\oplus\left( \bigcup_{i=1}^p a_i K_{n_i} \right),$$ for integers  $p, a_1, a_2,  \ldots, a_p, n_0, n_1,  \ldots, n_p \geq 1$ and $n_i \neq n_j$, for all $1\leq i < j \leq p$.
Thus $G$ has a total of $\sum_{i=1}^{p}a_i$ satellites   of  $p$ different orders.
In particular, the  graphs in Theorem \ref{theo:chordal2main} have at least two satellites. The ones of item (a) are  generalized core-satellite graphs with $p=2$ and $a_1=a_2=1$, while the graphs of item (b) are  core-satellite graphs having at least two satellites, that is,  $p=1$ and $a_1=t>1$. Both types of graphs have two main $\mathbf{Q}$-eigenvalues. We notice that a generalized core-satellite with one  satellite  is a complete graph. Hence, we may assume that either $p>1$, or that $p = 1$ and $a_1 \geq 2$.





In the sequence,  for a given generalized core-satellite graph $G$, we determine the exact number $k\geq 2$  of main $\mathbf{Q}$-eigenvalues of $G$.

 \begin{theorem}\label{theo:core-satkmain} Let $G$ be a generalized core-satellite graph.  Then:
\begin{enumerate}
    \item[\rm{(a)}] $G$ has  $k=2$ main $\mathbf{Q}$-eigenvalues if and only if $G$ has exactly two satellites, all with different orders, or $G$ has two or more satellites, all of the same order.

    \item[\rm{(b)}] Let $k\geq 3$. Then $G$ has $k$ main $\mathbf{Q}$-eigenvalues if and only if $G$ has at least three satellites, divided into $k-1$ different orders.
\end{enumerate}
\end{theorem}

\begin{proof}
Let $G$ be a  generalized core-satellite graph $G=K_{n_0} \oplus H$, where $$H= \bigcup_{i=1}^p a_i K_{n_i},$$ for integers  $p, a_1, a_2,  \ldots, a_p, n_0, n_1, n_2,  \ldots, n_p \geq 1$ and $n_i \neq n_j$, for all $1\leq i < j \leq p$.
\vspace{0.5em}

\noindent \textit{(a)} By Theorem \ref{theo:chordal2main},  $G$ has two main $\mathbf{Q}$-eigenvalues if and only if  $p=2$ and $a_1=a_2=1$\ or   $p=1$ \ and \ $a_1\geq 2$.
\vspace{0.5em}

\noindent \textit{(b)} Assume $p \geq 2$.
 Firstly, we note that the cograph $H=a_1 K_{n_1} \cup a_2 K_{n_2} \cup \ldots \cup a_p K_{n_p}$ has exactly $p$ main $\mathbf{Q}$-eigenvalues, since main eigenvalues are considered to be distinct. Therefore, its complement graph $\overline{H} = a_1 \overline{K_{n_1}} \oplus a_2\overline{K_{n_2}} \oplus \ldots \oplus a_p\overline{K_{n_p}}$ also has $p$  main $\mathbf{Q}$-eigenvalues,  which are all non-zero by Lemmas \ref{lem:zeroasQeigenval} and \ref{lem:cogbip},  since $\overline{H}$ is not bipartite for $p\geq 3$ or $p=2$ and $a_1\geq 2$ or $a_2 \geq 2$. Thus, $\overline{K_{n_0}}\cup \overline{H}$ has $p+1$ main $\mathbf{Q}$-eigenvalues, and so has its complement, which is  the initial graph $G=K_{n_0} \oplus H$. We conclude that  the generalized core-satellite graph $G$ has $k \geq 3$ main $\mathbf{Q}$-eigenvalues if and only  if  $k=p+1$, for $p= 2$ and $a_1 \geq 2 $ or $a_2 \geq 2$ or $p \geq 3$, that is, $G$ has at least three satellites of $k-1$ different orders.

\end{proof}

\section{Constructing \emph{quasi}-threshold non-generalized-core-satellite graphs with three main $\mathbf{Q}$-eigenvalues}\quad

In \citep{VINAGRE202033}, a recursive method to construct all threshold graphs with $k+1$ main  $\mathbf{Q}$-eigen\-values from threshold graphs having $k$ main $\mathbf{Q}$-eigenvalues was presented.  It is a natural question to consider whether a similar procedure to construct quasi-threshold graphs with $k+1$ main  $\mathbf{Q}$-eigenvalues from quasi-threshold graphs having $k$ main $\mathbf{Q}$-eigenvalues may be feasible.

In spite of the fact that we have a complete characterization of the basis case ($k=2$), as well as the characterization of all generalized core-satellite graphs -- a subclass of \emph{quasi}-threshold graphs --  having $k \geq 2$ main $\mathbf{Q}$-eigenvalues (see Theorem \ref{theo:core-satkmain}), which are steps towards the investigation of main $\mathbf{Q}$-eigenvalues of cographs, a similar construction does not seem possible for \emph{quasi}-threshold graphs in general. To justify this assertion, we present in the sequel several examples of \emph{quasi}-threshold graphs, which are not generalized core-satellite graphs, having  three main $\mathbf{Q}$-eigenvalues. We understand that these examples show that a complete description of all \emph{quasi}-threshold graphs having $k>2$ main $\mathbf{Q}$-eigenvalues is a very hard problem.

In order to present our constructions,  we  start by recalling that,  by Proposition~\ref{prop:qt-con},  a connected non-complete \emph{quasi}-threshold graph can be described as $G=K_c\oplus H$, for some integer $c \geq 1$ and a disconnected \emph{quasi}-threshold graph $H$. Also,  from Theorem~\ref{teo:qt-estrut}, we know that, if such graph $G$ has three main $\mathbf{Q}$-eigenvalues then the graph $H$ has two main $\mathbf{Q}$-eigenvalues, if we take Lemma \ref{lem:maincompl} into account.

To present our constructions, we recall some facts found in the literature that we state as the following remark for further reference.

\begin{remark}\label{rem:conditions}
  The $\mathbf{Q}$-spectrum of the split complete graph  $K_a \oplus \overline{K_b}$, where $a \geq 1$ and $b>1$   is
  \[\{(a+b-2)^{[a-1]},a^{[b-1]}\}\cup \{q_1^*, q_2^*\}\,\]
  where $q_1, q_2 $ are the two roots of the equation $q^2-(b+3a-2)q+(2a^2-2a)=0$. Since $a+b-2$ and $a$ are the non-main eigenvalues determined by Proposition \ref{prop:nonmaincogr},   the two main $\mathbf{Q}$-eigenvalues of the graph are $q_1$ and  $q_2 $ as indicated, considering that the graph is not regular.  For our purposes, we give below  the signless Laplacian condensed matrix of the graph, which is,
   \[\mathbf{C}(K_{a} \oplus \overline{K_b})=\left[
       \begin{array}{cc}
         2(a-1)+b & \sqrt{ab} \\
         \sqrt{ab} & a \\
       \end{array}
     \right]\,.
     \]
In \citep{DEFREITAS20102352}, the authors give   a condition on $a$ and $b$ under which the  eigenvalues $q_1$ and $q_2$ are integers and provide  some examples of such graphs, as the infinite family where $a=2s-1$ and $b=3s$,  for each integer $s \geq 1$:  indeed,   it holds that $q_1=8s-4$ and $q_2=s-1$, for $s \geq 1$ .

On the other hand,  from Lemma \ref{lem:mainUnion} we have that  the two main $\mathbf{Q}$-eigenvalues of the graph $K_c \oplus  (K_{a_1}\cup K_{a_2})$, where $c, a_1, a_2$ are positive integers with $a_1 \neq a_2$, are the roots
$q_1, q_2$  of a certain quadratic equation.  Also in this case, we find in \citep{DEFREITAS20102352} conditions on parameters $a_1, a_2, c$ which allow  the spectrum of this kind of graph to consist entirely of integers. For example, for each positive integer $s$,  if $c=s, \ a_1=s+1 $ and $a_2=s+2$ then $q_1=5s+2$ and $q_2= 2s$;  and  if $c=a_1=a_2=s$ then $q_1=5s-2$ and $q_2=2s-2$.
   \end{remark}

In the following, we apply the facts listed in Remark \ref{rem:conditions} to  describe  some disconnected \emph{quasi}-threshold graphs  with two  main $\mathbf{Q}$-eigenvalues,  satisfying  the conditions of item (a) of Theorem \ref{teo:qt-estrut}, that could replace graph $H$ in decomposition given in Proposition \ref{prop:qt-con}. Actually, in each item of the example below, we give an infinite family of such graphs. By using  these graphs $H$,  connected \emph{quasi}-threshold graphs of type $G=K_c \oplus H$  with  three main $\mathbf{Q}$-eigenvalues, which are not generalized core-satellite graphs, can be constructed.

\begin{example}\quad

\begin{enumerate}

\item[(a)] $ \mbox{ For integers }  \  p,\  a \geq 1\ \mbox{ and } \ \ b \geq 2\,,$  the graph
    \[ H_1(a,b,p)=\overline{K_a}  \cup pK_b\,,\] \   has  the main $\mathbf{Q}$-eigenvalues $2b-2$ and  $0$;
 \vskip.3cm

\item[(b)] $\mbox{ For integers }  \  p\geq 2 \mbox{ and }  a, b \geq 2\,,$
    \[\label{eq1.2} H_2(a,b,p)= p(K_a\oplus \overline{K_b}) \]
has  two main $\mathbf{Q}$-eigenvalues in accordance to   Remark \ref{rem:conditions}.
  In the particular case
\[ H_2^{\prime}(b,p)=H_2(1,b,p)=  p(K_1\oplus \overline{K_b})\,, \]
for integers \ $ p\geq 2 \mbox{ and }  \ b \geq 2$,   we have the union of at least two bipartite chordal cographs with $1+b$ and 0 as their main $\mathbf{Q}$-eigenvalues. And from this latter case, we also  obtain
\[ H_2^{\prime \prime}(b,p_1,p_2)=\overline{K_{p_1}} \cup p_2(K_1\oplus \overline{K_b})\,, \] for positive  integers $  p_1 \mbox{ or } p_2 \geq 2 \  \mbox{ and } \ \ b  \geq 2\,$, with  the same main $\mathbf{Q}$-eigenvalues.
 \vskip.3cm

\item[(c)] $\ \mbox{ For integers }  \  p  \geq 2\,, \    s \geq 1 \  \ \mbox{ and } \ \  a_1\geq 2 \  \mbox{ or}  \ \ a_2 \geq 2$, the graph
    \[ H_3(s,a_1, a_2,p)= p(K_{s}\oplus \left(K_{a_1}\cup K_{a_2})\right) \, \  \]
has   two main $\mathbf{Q}$-eigenvalues as we have seen   in Lemma  \ref{lem:mainUnion};
 \vskip.3cm


\item[(d)] By choosing appropriate values for the parameter $b$ in the above
graph of type $H_2^{\prime}(b,p)$, we obtain, for example, $\  \ \mbox{ for positive  integers }  \ a \geq 2, p_1,  p_2 $, the graph  
\[ H_4(a,p_1,p_2)= p_1K_{a} \cup p_2(K_1\oplus \overline{K_{(2a-3)}})\,, \]
for which the main $\mathbf{Q}$-eigenvalues are $2a-2$ and 0, by Lemma 8.
 \vskip.3cm

\item[(e)]  From the two previous cases,  
     for positive  integers  $  a \geq 2, p_1, p_2, p_3$,
    we can define
    \[ H_5(a,p_1,p_2,p_3)= p_1K_{a} \cup p_2(K_1\oplus \overline{K_{(2a-3)}})\cup \overline{K_{p_3}}\,, \  \ \]  for which the  main $\mathbf{Q}$-eigenvalues are also $2a-2$ and 0;
\vskip0.3cm

\item[(f)] By taking into account    the conditions presented in Remark \ref{rem:conditions}, we can choose appropriate values for the parameters, to obtain, for example, the graph \[  H_6(s,p_1, p_2, p_3)= p_1(K_{(2s-1)}\oplus \overline{K_{(3s)}})\cup p_2K_{(4s-1)} \cup p_3K_{(s+1)/2}\,, \] where $  \  s,  p_1, p_2, p_3$ are positive integers and  $s$  is odd, which has $8s-4$ and $s-1$ as  main $\mathbf{Q}$-eigenvalues;
\vskip0.3cm

\item[(g)] Following the same ideas, we obtain   \[ H_7(s,p_1,p_2,p_3)=  p_1(K_{s}\oplus (K_{s+1}\cup K_{s+2})) \cup p_2K_{(5s+4)/2}\cup  p_3K_{(s+1)}\,,\]  $ \mbox{ for integers }  \  s, p_1, p_2, p_3  \geq 1 \ \mbox{ such that } \ \   s  \mbox{ is even}\,$,  for which the  main $\mathbf{Q}$-eigenvalues are $5s+2$ and $2s$, by Remark \ref{rem:conditions};
\vskip0.3cm

\item[(h)] Still from the conditions presented in Remark \ref{rem:conditions}, for integers $  s, p_1, p_2, p_3  \geq 1$\   such that $s$ is even,  we obtain the graph
\[ H_8(s, p_1,p_2,p_3)= p_1(K_{s}\oplus (K_{s}\cup K_s)) \cup p_2K_{(5s/2)}\cup  p_3K_{s}\,, \]
for which the  main $\mathbf{Q}$-eigenvalues are $5s-2$ and $2s-2$.

\end{enumerate}
\end{example}

With the above example, we hope to have shown that a general characterization of \emph{quasi}-threshold graphs with $k\geq 3$ main $\mathbf{Q}$-eigenvalues, in the way it was done for threshold graphs in \citep{VINAGRE202033} and above, in the case of generalized core-satellites,  does not seem to be possible.
Even the case of connected cographs with $k=2$ main $\mathbf{Q}$-eigenvalues remains open. We note that the structural characterization of \emph{quasi}-threshold graphs presented in Proposition \ref{prop:qt-con}, which was based on the description of the minimal separators of vertices presented in Lemma \ref{lem:abreuchordal} of \citep{ABREU202168}, was fundamental for characterizing \emph{quasi}-threshold connected graphs in Theorem \ref{theo:chordal2main}. Therefore, in a possible approach to the case of connected cographs with $k=2$ main $\mathbf{Q}$-eigenvalues, it seems important to obtain a similar structural description.


\bibliographystyle{unsrtnat}
\bibliography{biblio2023}  






\end{document}